\documentclass{article}
\usepackage{amssymb}
\usepackage{amsmath}
\usepackage{amsthm}
\usepackage{mathrsfs}
\usepackage{mathtools}
\usepackage[hidelinks]{hyperref}
\usepackage[a4paper,hmargin=3cm]{geometry}
\IfFileExists{newcommand.tex}{\newcommand{\EC}{\mathbb{C}}

\newcommand{\MB}{\mathcal{B}}

\newcommand{\ijbar}{{i\overline{j}}}

\newcommand{\MQ}{\mathcal{Q}}

\newcommand{\MG}{\mathcal{G}}

\newcommand{\trace}{{\text{tr}}}

\newcommand{\ER}{\mathbb{R}}
\newcommand{\ED}{\mathbb{D}}

\newtheorem{corollary}{Corollary}[section]}{%
  \newcommand{\EC}{\mathbb{C}}
  \newcommand{\MB}{\mathcal{B}}
  \newcommand{\ijbar}{{i\overline{j}}}
  \newcommand{\MQ}{\mathcal{Q}}
  \newcommand{\MG}{\mathcal{G}}
  \newcommand{\trace}{{\text{tr}}}
  \newcommand{\ER}{\mathbb{R}}
}

\numberwithin{equation}{section}
\newtheorem{theorem}{Theorem}[section]
\newtheorem{lemma}[theorem]{Lemma}
\newtheorem{prop}[theorem]{Proposition}
\mathtoolsset{showonlyrefs}
\title{An Optimal Interior $C^{2,\beta}$ Estimate for  the Complex Monge-Amp\`ere Equation}
\author{Jie Deng,\  Jingchen Hu,\ Liding Huang, \ Jiawei Liu,\  Li Sheng,\  Xi Zhang}

\begin{document}
\maketitle

\begin{abstract}
We prove an interior ${}C^{2, \beta}$ estimate for plurisubharmonic solutions
$u\in C^{1,\alpha}$ of the complex Monge-Amp\`ere equation
$\det(u_{\ijbar})=f$ when  ${}f>0$, ${}f\in C^{\beta}$ and \({}\alpha>1-2/n\), optimizing the ${}C^{2,\beta}$ estimates of \cite{DZZ2011} and \cite{LLZ2018}.

When the righthand side  $f$ fails to be $C^{\beta}$, but the Euclidean Laplacian of \(f^{1/n}\)
is  bounded from below in the sense of distributions, we can prove a ${}C^{1,1} $ interior estimate.

For a smooth positive \({}f\), the estimate implies that every such solution
is  smooth, proving the conjecture of B{\l}ocki and Dinew
posed in~\cite[Question~24]{DGZ2016}.
\end{abstract}

\section{Introduction}
A plurisubharmonic solution of the complex Monge-Amp\`ere equation
\begin{equation}
\det(u_{\ijbar})=f\label{complexMongeAmpere}
\end{equation}
can
have an unbounded Hessian even when $f$ is smooth and bounded away
from zero.  Wang and Wang even gave an example on ${}\ED\times \EC^{n-1}$ that fails to be Dini-continuous \cite{WangWang2025}.

B{\l}ocki and  Dinew asked whether $u\in C^{1,\alpha}$ with exponent \(\alpha>1-2/n\) forces interior smoothness \cite[Question 24]{DGZ2016}. Actually 
Dinew and Dinew proved  that if 
$\det(u_{\ijbar} ) \geq  1$, and $u\in C^{1,\alpha}$ with $\alpha> 1-\frac{2}{n}$, then the zero set of $u$ contains no positive-
dimensional analytic subset(\cite{DinewDinew2016}, Theorem 35). They suggested that the same gradient hypothesis should imply
classical regularity for positive H\"older densities (\cite{DinewDinew2016}, Theorem 2 and the discussion following
it).

The threshold ${}\alpha>1-2/n$  is optimal.
Indeed, with \({}z'=(z_2,\ldots,z_n)\), the function \({}(1+|z_1|^2)|z'|^{2-2/n}\)
belongs to \(C^{1,1-2/n}\)  for  ${}n\geq 3$ and has smooth positive right-hand side ${}f$, but its Hessian
is unbounded near ${}z'=0$~\cite[p.~412]{BD2011}. In dimension two the same function is Lipschitz and fails to be differentiable  on \({}z'=0\).

In 2011 Dinew, Zhang and Zhang \cite{DZZ2011} proved  that
${}C^{1,1}$ solutions with positive \({}f\in C^\beta\)  are locally
$C^{2,\beta}$ for \(0<\beta<1\).
In 2018 Li, Li and Zhang \cite[Theorem~5.1(2)]{LLZ2018} obtained an interior
$C^{1,1}$ estimate from $u\in C^{1,\alpha}$, \({}\alpha>1-1/n\),
for \({}f\in C^{2,1}\) bounded  away from zero.
For  positive \({}f\in C^\beta\), \(0<\beta<1\), the interior ${}C^{2,\beta}$
estimate requires \({}1-\beta/[n(2+\beta)-1]<\alpha<1\).
Recently the global ${}C^{2,\beta}$ estimate for complex Monge-Amp\`ere equation with the righthand side ${}f\in C^{\beta}$ is proved over compact manifold \cite{BingWang26} and in a bounded strongly pseudoconvex domain with $C^3$ boundary \cite{ZhangXi-HuangliDing}.

For the real Monge-Amp\`ere equation,
Urbas~\cite{U1988} proved the corresponding result: a convex \({}C^{1,\alpha}_{}\) Aleksandrov solution of $\det D^2u=f>0$ with smooth ${}f$ is smooth in the interior if ${}\alpha>1-2/n$.

In this paper, since we only study interior  regularity, it suffices to work in ${}B_2\subset \EC^n$. 
We denote the full $C^{k,\sigma}(\Omega)$ norm by ${}|w|_{k,\sigma; \Omega}$,
writing ${}|w|_{\sigma; \Omega}$ when ${}k=0$. We denote the H\"older
seminorm of $w$ by $[w]_{\sigma; \Omega}$. In particular,
\[
 {}|w|_{1,1; \Omega}
 =\sup_{\Omega}|w|+\sup_{\Omega}|Dw|+[Dw]_{1; \Omega},
 \qquad
 [Dw]_{1; \Omega}
 =\sup_{\substack{x,y\in\Omega\\x\ne y}}\frac{|Dw(x)-Dw(y)|}{|x-y|}.
\]
For H\"older densities the two-scale comparison in Section~\ref{sec:holderDensity}
gives the following interior estimate.
\begin{theorem}[Hölder Right-Hand Side]\label{thm:holderDensity}
	Let $n\ge 2$, $1-2/n<\alpha\le 1$, and ${}0<\beta<1$.
	Suppose ${}u\in C^{1,\alpha}(B_2)$ is plurisubharmonic
	and solves \eqref{complexMongeAmpere} in the Bedford-Taylor sense.
	If ${}f\in C^\beta(B_2)$ and ${}f_\ast:=\inf_{B_2} f>0$, then ${}u\in C^{2,\beta}(B_1)$ and
	\[
	{}|u|_{2,\beta; B_1}\le C.
	\]
	The constant depends  only on ${}n,\alpha,\beta$, $|u|_{1,\alpha; B_2}$,
	${}|f|_{\beta; B_2}$, and ${}f_\ast$.
\end{theorem}
When the righthand side ${}f$ fails to be H\"older continuous, we can relay on the Laplacian lower bound ${}f^{1/n}$ and get the following ${}C^{1,1}$ interior estimate.
\begin{theorem}[Interior ${}C^{1,1}$ estimate]\label{thm:hessian}
Let $B_2=B_2(0)$, \({}n\ge2\),  and \(1-2/n<\alpha \leq 1\). Suppose
\(u\in C^{1,\alpha}(B_2)\) is plurisubharmonic and solves \eqref{complexMongeAmpere}
in the Bedford-Taylor sense. Assume that \(f\in L^1(B_2)\). Denote that ${}g=f^{1/n}$.
Assume that for some ${}g_*>0$ and \({}L\ge0\),
\begin{equation}\label{densityLaplacianAssumption}
       g\ge g_*\quad\text{a.e.},\qquad \Delta_e g\ge-L
       \quad\text{in distributions on }B_2.
\end{equation}
Then ${}u\in C^{1,1}(B_1)$ and
\begin{equation}\label{realHessianEstimate}
 |u|_{1,1; B_1}\le C.
\end{equation}
The  constant depends only on \({}n,\alpha,|u|_{1,\alpha; B_2},
g_*,L\), and the mean of \(g\) on ${}B_2$.
\end{theorem}

As an application of the estimates above, 
we give an affirmative answer to the question of  B{\l}ocki and Dinew recorded
in Dinew, Guedj and Zeriahi's 
survey \emph{Open problems in
pluripotential theory}~\cite[Question~24]{DGZ2016}.

\begin{corollary}[Smooth Right-Hand Side]\label{thm:main}
Let ${}n\geq2$ and ${}1-2/n<\alpha\le1$. Suppose  ${}u\in C^{1,\alpha}(\Omega)$ is plurisubharmonic and solves  \eqref{complexMongeAmpere} with $f\in C^\infty(\Omega)$, \(f>0\). Then $u\in C^\infty(\Omega)$.
\end{corollary}
The idea in  the proof of Theorem \ref{thm:holderDensity} and \ref{thm:hessian} are very similar. We let a molification  of the solution $u$: ${}C_3(\MG_tu+A\MG_hu-C_1-C_2|z|^2)$ touches ${}u$ from below. The distance of this $\MG_tu+A\MG_hu$ and $u$ contains  the information of the Euclidean Laplacian of  $u$. The contact point turns out to be roughtly the maximum  point of Laplacian ${}u$. So controling the Hessian lower bound of ${}\MG_tu+A\MG_hu$ gives a control of Hessian of ${}u$.

We first introduce the following Green  averaging operator (also see \cite{HS2026}).
\section{The Green average and its complex Hessian}
Since  we are only interested about local regularity in this paper, we study solutions in ${}B_2$.
The estimate in $B_{2r}(p)$ follows by a rescaling. We write
${}\Delta_e=4\sum_j\partial_j\partial_{\bar j}$.
Under \eqref{densityLaplacianAssumption}, use the upper semicontinuous representative of \(g\) obtained by subtracting
\(L|z|^2/(4n)\) from the subharmonic representative of \({}g+L|z|^2/(4n)\).
Its small-ball means converge  to its value at every point, so $g\ge g_*$
everywhere. For $z\in B_1$, the submean inequality on \(B_{1/2}(z)\)
bounds ${}g(z)$ by $4^{2n}$ times the mean of \(g\) on \({}B_2\), plus
${}C_n L$.

For a continuous function $w$ on ${}B_2$, write \({}\MB_s w(z)\) and
${}\mathcal{S}_s w(z)$ for its Euclidean  averages over ${}B_s(z)$ and
$\partial B_s(z)$, respectively. For $z\in\overline{B_1}$ and ${}0<t\le1/2$, set
\begin{equation}\label{greenAverageDefinition}
 {}\MG_t w(z)=\frac{2}{t^2}\int_0^t s\MB_s w(z)\,ds,\qquad
 \MQ_t w(z)=\MG_t w(z)-w(z).
\end{equation}
We use the same averages for locally integrable functions whenever the
integrals exist. For the  solution $u$ write ${}H_t=\big((\MG_t u)_{\ijbar}\big)$.
Throughout the  proof \({}\delta=1-\alpha\) and \(g=f^{1/n}\).
All bounds on ${}B_1$ are uniform for \({}0<t\le1/2\).

\begin{lemma}\label{lem:green}
Under the hypotheses of Theorem~\ref{thm:hessian}, ${}\MG_t u$ is $C^2$ and plurisubharmonic on $B_1$. There are
constants \(0<\kappa\le1\le K\), depending only on $n$, $\alpha$,
${}[Du]_{\alpha; B_2}$, ${}g_*$, ${}L$, and the mean of $g$ on ${}B_2$, such that
\begin{align}
 0\le \MQ_t u\le K t^{2-\delta},\qquad
 \kappa t^{(n-1)\delta} I\le H_t\le K t^{-\delta}I, \label{averageErrorAndHessianBounds}\\
 \det H_t\ge (\MG_t g)^n\ge f-K t^2,\qquad \det H_t\ge \kappa. \label{averagedDeterminantLowerBound}
\end{align}
Moreover,
\begin{equation}\label{sphericalMeanIdentity}
 \MQ_t u\le \mathcal{S}_tu-u=\frac{t^2}{n}\trace\,H_t.
\end{equation}
\end{lemma}

\begin{proof}
Subharmonicity gives $u\le \MG_t u\le \mathcal{S}_tu$, since the spherical mean is
nondecreasing in the radius. The affine Taylor polynomial at \({}z\) has average
$u(z)$, and its remainder on \(B_t(z)\) is ${}O(t^{1+\alpha})$. Thus \(\MQ_t u\le K t^{2-\delta}\) for large enough \({}K\).

For real coordinate indices \(i,j\), the divergence theorem gives
\begin{equation}\label{ballAverageDerivative}
 {}D_iD_j\MB_s u(z)=\frac{1}{|B_s|}
 \int_{\partial B_s(0)}[D_j u(z+h)-D_j u(z)]\,\mathbf n_i(h)\,dS(h),
\end{equation}
where \(\mathbf n\) is the outer  unit normal and $|B_s|$ is Euclidean
volume. This  expression is continuous and is $O(s^{-\delta})$.
Integration in \(s\) gives ${}\MG_t u\in C^2$ and, after increasing $K$, \({}H_t\le K t^{-\delta}I\).
Averages of  translates preserve plurisubharmonicity.

Fix a constant Hermitian matrix $M>0$ with \(\det M=1\). The mixed
Monge-Amp\`ere inequality~\cite{D2009} for $u$ and a quadratic potential
with complex Hessian $M^{-1}$ gives ${}\trace(M(u_{\ijbar}))\ge ng$ in distributions.
Averaging  gives $\trace(MH_t)\ge n\MG_t g$ for every such \(M\). Hence, pointwise,
\begin{equation}\label{determinantTraceInfimum}
 {}(\det H_t)^{1/n}
 =\inf_{\substack{M=M^*>0\\ \det M=1}}\frac{1}{n}\trace(MH_t)
 \ge \MG_t g .
\end{equation}

At each center ${}z$, the function \(g(w)+L|w-z|^2/(4n)\) is subharmonic.
Its ball submean inequality and $\MB_s(|\,\cdot-z|^2)(z)=ns^2/(n+1)$ give
\begin{equation}\label{densityAverageLowerBound}
 \MB_s g(z)\ge g(z)-\frac{Ls^2}{4(n+1)},\qquad
 \MG_t g(z)\ge g(z)-\frac{Lt^2}{8(n+1)}.
\end{equation}
Since $g$ is bounded on \({}B_1\), the latter inequality gives
\((\MG_t g)^n\ge f-Kt^2\) after increasing ${}K$. Also  \(\MG_t g\ge g_*\),
so $\det H_t\ge g_*^n$. The upper eigenvalue bound  in \eqref{averageErrorAndHessianBounds}
now gives the lower bound there with $\kappa=\min\{1,g_*^n\}/K^{n-1}$.

For smooth $u$, the divergence theorem  gives
$\frac{d}{ds}\mathcal{S}_su=\frac{s}{2n}\MB_s(\Delta_e u)$, where \(\Delta_e\)
is the real Euclidean Laplacian. Integration yields
$\mathcal{S}_tu-u=t^2\Delta_e \MG_t u/(4n)$, as  in~\cite[(3.3)-(3.4)]{HS2026}.
Local mollification gives the same identity for \({}u\), since $\MG_t u$
is \(C^2\) at each fixed \({}t\). Now \(\Delta_e \MG_t u=4\trace\,H_t\).
\end{proof}
\section{Interior \texorpdfstring{$C^{2,\beta}$}{C2,beta} estimate for a H\"older density}
\label{sec:holderDensity}
In the  estimates below \({}C\) depends
only on \(n,\alpha,\beta,|u|_{1,\alpha; B_2},
|f|_{\beta; B_2}\), and \(f_\ast\). A constant with subscript \(\theta\) also
depends on \(\theta\).

\subsection{An estimate for the averaging error}

\begin{prop}\label{prop:holderContact}
	Under the hypotheses of Theorem~\ref{thm:holderDensity}, suppose \(0<\delta=1-\alpha<2/n\). Set
	\[
	{}\tau=1-\frac{\beta}{2},\qquad
	\gamma=\frac{\beta\delta}{\tau(2-n\delta)} .
	\]
	For \(0<\theta<1\), \({}B_r(p)\subseteq B_1\), \(0<r\le1\), and
	\(0<t<1/2\), one has
	\begin{equation}\label{holderLocalizedContact}
		\MQ_t u(p)\le \theta\sup_{B_r(p)}\MQ_t u
		+ C_\theta r^{-\gamma}t^{2-\tau\delta}.
	\end{equation}
\end{prop}

\begin{proof}
	For \(0<s\le1/2\) and centers  in \(\overline{B_1}\), the Taylor,
	determinant, and trace estimates in the proof of Lemma~\ref{lem:green} give
	\begin{align*}
		0\le\MQ_s u&\le Ks^{2-\delta},&
		\kappa s^{(n-1)\delta} I&\le H_s\le Ks^{-\delta}I,\\
		\det H_s&\ge f_\ast,&
		\MQ_s u&\le \frac{s^2}{n}\trace H_s .
	\end{align*}
	We may take
	\[
	K=C_{n,\alpha}(1+[Du]_{\alpha; B_2})\ge1,\qquad
	\kappa=\min\{1,f_\ast\}K^{1-n}.
	\]
	Since
	\({}[g]_{\beta; B_2}\le f_\ast^{1/n-1}[f]_{\beta; B_2}/n\),
	\eqref{determinantTraceInfimum} gives
	\[
	{}(\det H_s)^{1/n}\ge\MG_s g\ge g-[g]_{\beta; B_2}s^\beta,
	\qquad \det H_s\ge f-Cs^\beta .
	\]
	
	Put \(h=t^\tau\), \(M=\sup_{B_r(p)}\MQ_t u\), and
	\({}c_\theta=\min\{1/4,\theta\kappa/(4K)\}\).
	If \(M\le Kt^{2-\tau\delta}\),  then \eqref{holderLocalizedContact} follows
	from \(r\le1\). If \({}h^{2-n\delta}>c_\theta r^2\), then
	\(t^{-\beta\delta/2}\le c_\theta^{-\gamma/2}r^{-\gamma}\). In this case
	\[
	{}\MQ_t u(p)\le Kt^{2-\delta}
	=Kt^{2-\tau\delta}t^{-\beta\delta/2}
	\le C_\theta r^{-\gamma}t^{2-\tau\delta}.
	\]
	
	It remains to consider \({}M>Kt^2h^{-\delta}\) and
	\(h^{2-n\delta}\le c_\theta r^2\),  which gives \({}t<h\le1/2\).
	Set \(A=\theta M/(Kh^{2-\delta})\). Since \(t^2/h^2=t^\beta\)
	and \({}M\le Kt^{2-\delta}\), we have
	\[
	{}\theta t^\beta<A\le\theta(t/h)^{2-\delta}\le\theta,\qquad
	A\MQ_h u\le\theta M,\qquad
	AH_h\ge\frac{\theta\kappa M}{K h^{2-n\delta}}I\ge\frac{4M}{r^2}I .
	\]
	The function
	\({}F(z)=\MQ_t u(z)+A\MQ_h u(z)-2Mr^{-2}|z-p|^2\)
	is nonnegative at \({}p\) and is  at most \((\theta-1)M\) on
	\({}\partial B_r(p)\). At an interior maximum \({}\xi\),
	\[
	{}\MQ_t u(p)\le F(\xi)\le \MQ_t u(\xi)+\theta M.
	\]
	As in Proposition~\ref{prop:contact},  the function
	\[
	{}\phi(z)=\frac{\MG_t u(z)+A\MG_h u(z)-2Mr^{-2}|z-p|^2-F(\xi)}{1+A}
	\]
	touches \(u\) from below at \(\xi\). At this point
	\((1+A)(\phi_{\ijbar})\ge H_t+(A/2)H_h>0\), so comparison gives
	\[
	\det\bigl({}H_t+(A/2)H_h\bigr)\le (1+A)^n f(\xi).
	\]
	
	Evaluate the matrices at \(\xi\) and put \({}R=H_t^{-1/2}H_hH_t^{-1/2}\).
	Since \(\det(I+(A/2)R)\ge1+(A/2)\trace R\),  we obtain
	\[
	\frac A2\det H_t\,\trace R
	\le (1+A)^n f(\xi)-\det H_t
	\le C(A+t^\beta).
	\]
	The lower bound \({}A>\theta t^\beta\) implies
	\(\det H_t\,\trace R\le C_\theta\). Applying the arithmetic-geometric
	mean inequality to the \(n-1\) largest eigenvalues of \(R\) gives
	\[
	\lambda_{\min}(R)
	\ge \frac{(n-1)^{n-1}\det R}{(\trace R)^{n-1}}
	=\frac{(n-1)^{n-1}\det H_h(\det H_t)^{n-2}}
	{(\det H_t\,\trace R)^{n-1}}
	\ge C_\theta^{-1}.
	\]
	Consequently \({}H_t\le C_\theta H_h\le C_\theta Kh^{-\delta}I\). The
	trace estimate yields \(\MQ_t u(\xi)\le C_\theta t^2h^{-\delta}
	=C_\theta t^{2-\tau\delta}\). The bound for \({}F(\xi)\) now gives
	\eqref{holderLocalizedContact}.
\end{proof}

Proposition \ref{prop:holderContact} implies the following estimate for ${}\MG_t$:
\begin{prop}
	Suppose that $u$ is a  ${}C^{1,1-\mu}$ solution to equation (\ref{complexMongeAmpere}) in $B_2$, under the hypotheses of Theorem~\ref{thm:holderDensity}, in particular $0<\mu<2/n$. Let ${}\tau=1-\beta/2$. Then in ${}B_{1/2}$, with $t\leq 1/2$, 
		\begin{equation}
	\MQ_t u\le C t^{2-\tau \mu}.\label{holderOneStepGain}
	\end{equation}
	The constant $C$ depends on
	\(n,\mu,\beta,f_\ast\), \({}|f|_{\beta; B_2}\), \(|u|_{1,1-\mu; B_2}\).
	\end{prop}
\begin{proof}
	Suppose that on the unit ball $B_1$, \({}(1-|z|)^\gamma\MQ_t u(z)\) achieves its maximum  at \({}p\),
	where $\gamma=\beta \mu/[\tau(2-n\mu)]$. This is the ${}\gamma$ in Proposition \ref{prop:holderContact} with ${}\delta$ replaced by ${}\mu$. Then apply  Proposition \ref{prop:holderContact} in $B_{2r}(p)$ with
	\({}r=(1-|p|)/2\)  and \({}\theta=2^{-\gamma-1}\) in
	\eqref{holderLocalizedContact}. Then
	\begin{align}
		(1-|p|)^{\gamma}\MQ_tu(p)\leq (1-|p|)^{\gamma} (\theta\max_{B_r(p)}\MQ_tu+C_\theta r^{-\gamma} t^{2-\tau\mu}).
	\end{align}
	If $\max_{B_r(p)}\MQ_tu$ is achieved at $q$, then 
	\begin{align}
		{}(1-|q|)^{\gamma}\max_{B_r(p)}\MQ_tu\leq 	(1-|p|)^{\gamma}\MQ_tu(p).
	\end{align} Combining the two inequalities above and use
	\begin{align}
		\frac{{}(1-|p|)^\gamma}{(1-|q|)^\gamma}\leq 2^\gamma.
	\end{align} We get 
	\begin{align}
		{}(1-|p|)^\gamma \MQ_tu(p)\leq 2^\gamma (1-|q|)^\gamma \MQ_tu(q)\theta+2^\gamma C_\theta t^{2-\tau\mu}\leq 2^\gamma (1-|p|)^\gamma \MQ_tu(p)\theta+2^\gamma C_\theta t^{2-\tau\mu}.
	\end{align}
	By choosing ${}\theta=2^{-\gamma-1}$, we get that $(1-|p|)^\gamma \MQ_tu(p)$ is bounded by \(Ct^{2-\tau \mu}\). It follows that ${}\MQ_tu$ in ${}B_{1/2}\leq Ct^{2-\tau\mu}$.
	\end{proof}
	For \({}0<a<1\), a local bound \(\MQ_t u\le Ct^{1+a}\) implies
	\(u\in C^{1,a}_{\mathrm{loc}}\). 
	This follows from Lemma~\ref{lem:greenErrorGradient} in the appendix.
	
	Thus one application of the previous Proposition  improves ${}u\in C^{1,1-\mu}(B_2)$ to $u\in C^{1,1-\mu\tau}(B_{1/4})$. To $u\in C^{1,1-\delta}$, we apply the estimate finite times on concentric balls and get that for  a very small ${}r_0$, ${}u\in C^{1,1-\tau^k\delta}(B_{r_0})$, with \(1-\beta/[n(2+\beta)-1]<1-\tau^k\delta\). Then we can apply Li, Li and Zhang's interior estimate
	\cite[Theorem~1.2]{LLZ2018} and get  \(|u|_{2,\beta; B_{r_0/2}}\le C\).
	
	Actually we can show that $u$ is $ C^{2,\beta}$ at any interior point of $B_2$, since using a biholomorphic automorphism, we can map any interior point of $B_2$ to ${}0$. In particular, $|u|_{2,\beta;B_1}$ has a uniform bound.

\begin{proof}[Proof of Corollary~\ref{thm:main}]
On every ball compactly contained in \({}\Omega\), the smooth positive density
\(f\) belongs to \({}C^{1/2}\) and is bounded away from zero.
Applying Theorem~\ref{thm:holderDensity} on such balls gives
\(u\in C^{2,1/2}_{\mathrm{loc}}(\Omega)\). Hence
\(u\in W^{2,p}_{\mathrm{loc}}(\Omega)\) for every finite \(p>1\).
Taking \({}p>n(n-1)\), the regularity criterion of B{\l}ocki and Dinew
\cite[p.~412, (2)]{BD2011} yields \(u\in C^\infty(\Omega)\).
\end{proof}

\section{Interior \texorpdfstring{$C^{1,1}$}{C1,1} estimate}

\subsection{A bound for the averaging error}
The following proposition is an analog of Proposition \ref{prop:holderContact}. The idea of proof is also very similar. One difference is that in the following proof of Proposition \ref{prop:contact} the second averaging radius ${}h$ is fixed, while in the previous proof, ${}h$ is an explicit function of $t$. Another difference is that in the previous proof we use $g\in C^\beta$ and get $\MG_t g\geq g-C_2t^\beta$, while in the following proof of Proposition \ref{prop:contact}, we use the Laplacian lower bound of $g$ to obtain a lower bound for ${}\MG_t g$: $\MG_tg\geq g-C_1t^2$,  which also allows us to derive the desired estimate.
\begin{prop}\label{prop:contact}
Under  the hypotheses of Theorem~\ref{thm:hessian}, for every $0<\theta<1$
there is a constant ${}C_\theta$, depending only on ${}\theta$ and  \({}n,\alpha,|u|_{1,\alpha; B_2},
g_*,L\), and the mean of \(g\) on ${}B_2$
such that for $0<t<1/2$,
\begin{equation}\label{unitBallAveragingEstimate}
 \MQ_t u(0)\le\theta\sup_{B_1}\MQ_t u+C_\theta t^2 .
\end{equation}
\end{prop}

\begin{proof}
All averages are computed with centers in \(\overline{B_1}\) and have radius at most ${}1/2$. 
By the proof of Lemma~\ref{lem:green}, we may take
${}K=C_{n,\alpha}(1+[Du]_{\alpha; B_2})\ge1$ so that
${}\MQ_s u\le Ks^{2-\delta}$ and $H_s\le Ks^{-\delta}I$.
Since \(\det H_s\ge g_*^n\), we also have
\({}H_s\ge\kappa s^{(n-1)\delta}I\) with
$\kappa=\min\{1,g_*^n\}K^{1-n}$.
We  fix the radius $h$:
\begin{equation}\label{auxiliaryRadius}
 {}h=\left(\frac{\theta\kappa}{4K}\right)^{1/(2-n\delta)}.
\end{equation}
Since \({}0<2-n\delta\le2\) and \(\kappa\le1\le K\), we have ${}h\le\sqrt{\theta}/2<1/2$. By construction, $4Kh^{2-\delta}=\theta\kappa h^{(n-1)\delta}$.
If \({}h\le t<1/2\), the rough bound gives \({}\MQ_t u\le Kh^{-\delta}t^2\).

Suppose $0<t<h$ and put ${}M=\sup_{B_1}\MQ_t u$.
If ${}M\le t^2$, the conclusion is immediate. Otherwise set \({}A=4M/(\kappa h^{(n-1)\delta})\). Both \(M\) and ${}\MQ_h u$ are at most \(Kh^{2-\delta}\), so
\begin{equation}\label{weightedErrorBound}
 {}A\le\theta,\qquad A\MQ_h u\le\theta M,\qquad AH_h\ge4M I .
\end{equation}
The function $F=\MQ_t u+A\MQ_h u-2M|z|^2$ satisfies
${}F\le A\MQ_h u-M|z|^2\leq (\theta-1)M<0$ on ${}\partial B_1$ and \(F(0)\ge0\).
At an interior maximum ${}\xi$ we therefore have
\begin{equation}\label{contactPointComparison}
 {}\MQ_t u(0)\le F(\xi)\le\MQ_t u(\xi)+\theta M .
\end{equation}
The function
\begin{equation}\label{lowerTouchingFunction}
 {}\phi(z)=\frac{\MG_t u(z)+A\mathcal{G}_h u(z)-2M|z|^2-F(\xi)}{1+A}
\end{equation}
touches ${}u$ from  below at \({}\xi\). By  \eqref{weightedErrorBound},
\({}(1+A)(\phi_{\ijbar}(\xi))\geq H_t(\xi)+2M I>0\).
Bedford-Taylor comparison gives \(\det({}\phi_{\ijbar}(\xi))\le g(\xi)^n\). Indeed, if the inequality
failed, upper semicontinuity of \(g\) and continuity of the test Hessian would
give a ball ${}B_b(\xi)\Subset B_1$ on which
\(\phi+e(b^2/2-|z-\xi|^2)\) is strictly plurisubharmonic with determinant
greater than $g^n$ for some \(e>0\). This function is  below ${}u$ on the
boundary and above it at \({}\xi\), contrary to comparison~\cite{BT1982}.

All matrices below are evaluated at ${}\xi$. Taking the trace against
$(\phi_{\ijbar})^{-1}$ gives
\begin{equation}\label{contactTraceInequality}
 \trace\bigl((\phi_{\ijbar})^{-1}H_t\bigr)
 +2M\trace\bigl((\phi_{\ijbar})^{-1 }\bigr)\le n(1+A).
\end{equation}
By the arithmetic-geometric mean inequality and \eqref{densityAverageLowerBound},
\begin{equation}\label{contactTraceLowerBound}
 \trace\bigl((\phi_{\ijbar})^{-1}H_t\bigr)
 \ge n\left(\frac{\det H_t}{\det(\phi_{\ijbar})}\right)^{1/n}
 \ge n\frac{\MG_t g(\xi)}{g(\xi)}\ge n-Ct^2 .
\end{equation}
Subtracting gives
\begin{equation}\label{inverseHessianTraceBound}
 {}\trace\bigl((\phi_{\ijbar})^{-1}\bigr)\le\frac{nA+Ct^2}{2M}
 \le\frac{C}{\kappa h^{(n-1)\delta}},
\end{equation}
since \(A/(2M)=2/(\kappa h^{(n-1)\delta})\), ${}M>t^2$ and
$\kappa h^{(n-1)\delta}\le1$.
Each eigenvalue  of $(\phi_{\ijbar})$ is at most
$\det(\phi_{\ijbar})\bigl[\trace\bigl((\phi_{\ijbar})^{-1}\bigr)\bigr]^{n-1}$.
Hence  \({}(\phi_{\ijbar})^{-1}\ge c(\kappa h^{(n-1)\delta})^{n-1} I\).
By \eqref{sphericalMeanIdentity} and the upper trace bound,
\begin{equation}\label{contactErrorEstimate}
 \begin{aligned}
 \MQ_t u(\xi)&\le\frac{t^2}{n}\trace\,H_t(\xi)\\
 &\le C(\kappa h^{(n-1)\delta})^{1-n}t^2
     \trace\bigl((\phi_{\ijbar})^{-1}H_t\bigr)
 \le C(\kappa h^{(n-1)\delta})^{1-n}t^2 .
 \end{aligned}
\end{equation}
Together  with \eqref{contactPointComparison}, this proves
\eqref{unitBallAveragingEstimate}. The constants in the trace estimates depend
only on \({}n,g_*,L\) and $\sup_{{}B_1}g$.

For later use, the two ranges of ${}t$ give
$C_\theta\le C(1+Kh^{-\delta}+(\kappa h^{(n-1)\delta})^{1-n})$.
From \eqref{auxiliaryRadius} and $K\ge1$ we have
${}h\ge c_\theta K^{-n/(2-n\delta)}$.
Consequently
\begin{equation}\label{unitBallConstantDependence}
 C_\theta\le C_{\theta,0}\,K^{2(n-1)^2/(2-n\delta)},
\end{equation}
where ${}C_{\theta,0}$ depends  only on ${}n,\alpha,\theta,g_*,L$ and $\sup_{{}B_1}g$. 
\end{proof}

\begin{prop}[Rescaling to a ball of radius \(r\)]\label{prop:rescaledAveraging}
Under the hypotheses of  Theorem~\ref{thm:hessian}, set
\({}\gamma=2(n-1)^2\delta/(2-n\delta)\). For every \({}0<\theta<1\) there is
a constant \({}C_\theta\), depending only on $\theta$ and  \(n,\alpha,|u|_{1,\alpha; B_2},
g_*,L\), and  the mean of \({}g\) on ${}B_2$
such that for \({}\overline{B_r(p)}\subset B_1\), ${}0<r\le1$ and \(0<t<1/2\),
\begin{equation}\label{localSupremumEstimate}
 \MQ_t u(p)\le \theta\sup_{B_r(p)}\MQ_t u+C_\theta r^{-\gamma}t^2 .
\end{equation}
\end{prop}

\begin{proof}
Put ${}a=r/2$ and \(v(x)=a^{-2}u(p+ax)\) on \(B_2\). The rescaled n-th root of the right-hand side
\({}g(p+a\,\cdot)\) has the same upper and lower bounds and has Laplacian
at least $-{}a^2L$. Moreover, 
\begin{equation}\label{unitBallScaling}
 {}\MQ_s v(x)=a^{-2}\MQ_{as} u(p+ax),\qquad
 [Dv]_{\alpha; B_2}=a^{-\delta}[Du]_{\alpha; B_{2a}(p)}.
\end{equation}
If $t<a/2$, apply Proposition~\ref{prop:contact} with $s=t/a$ and multiply
by ${}a^2$. By \eqref{unitBallConstantDependence}, the rescaled constant is at
most ${}C_\theta a^{-\gamma}$. Since $a=r/2$ and \(B_a(p)\subset B_r(p)\),
this gives \eqref{localSupremumEstimate}. If ${}t\ge r/4$, then
$\MQ_t u\le Kt^{-\delta}t^2\le K(4/r)^\delta t^2$, which gives the
same estimate because $\gamma\ge\delta$ and \({}r\le1\).
\end{proof}

\begin{prop}[Interior Laplacian bound]\label{prop:interiorLaplacian}
Under the hypotheses of Theorem~\ref{thm:hessian}, there is a constant \({}C\)
depending only on ${}n$, ${}\alpha$, ${}[Du]_{\alpha; B_2}$, $g_*$, $L$,
and the mean of ${}g$ on ${}B_2$, such that for \(0<t<1/2\),
\begin{equation}\label{quadraticAveragingError}
 {}0\le \MQ_t u\le Ct^2 \qquad\text{on }B_1.
\end{equation}
Moreover, \({}0\le\Delta_e u\le 8(n+1)C\) in  distributions on $B_1$.
\end{prop}

\begin{proof}
When \(\delta=0\), Lemma~\ref{lem:green} gives \eqref{quadraticAveragingError}
directly. Suppose \(\delta>0\),  so $\gamma>0$. We first bound
\(\MQ_t u(0)\). Set \(\rho(z)=1-|z|\) on \({}\overline{B_1}\), which  is the distance to the boundary of ${}B_1$ and let
$p$ maximize $\rho^\gamma\MQ_t u$. If the maximum is positive,
then \(p\in B_1\). Since \({}\rho(z)\ge\rho(p)/2\) on  $B_{\rho(p)/2}(p)$,
\begin{equation}\label{weightedMaximumBound}
 \sup_{B_{\rho(p)/2}(p)}\MQ_t u\le 2^\gamma\MQ_t u(p).
\end{equation}
Apply \eqref{localSupremumEstimate} to this ball with ${}\theta=2^{-\gamma-1}$.
Moving ${}\MQ_t u(p)/2$ to the left gives
\({}\rho(p)^\gamma\MQ_t u(p)\le Ct^2\), which also holds if the maximum
is zero. As $\rho(0)=1$, this bounds ${}\MQ_t u(0)$.

For ${}x\in B_1$, apply the estimate at the origin to ${}v(y)=4u(x+y/2)$
on ${}B_2$. We  have $[Dv]_{\alpha; B_2}\le 2^\delta[Du]_{\alpha; B_2}$.
The function ${}g(x+\cdot/2)$ is bounded below by ${}g_*$, its Laplacian
is bounded below by $-{}L/4$, and its mean on ${}B_2$ is at most $2^{2n}$
times the mean of $g$ on ${}B_2$. Thus \(4\MQ_t u(x)=\MQ_{2t} v(0)\le C(2t)^2\)
for \(t<1/4\). For $1/4\le t<1/2$, Lemma~\ref{lem:green} gives
\({}\MQ_t u\le K4^\delta t^2\). This proves \eqref{quadraticAveragingError}.

For \(\psi\in C_c^\infty(B_1)\), the second-order Taylor formula gives
\({}(\MQ_t\psi)/t^2\to\Delta_e\psi/[8(n+1)]\) uniformly, since
the mean of each squared real coordinate on \({}B_s(0)\) is \({}s^2/(2n+2)\).
By symmetry of ${}\MG_t$,
\(\int\psi\,\MQ_t u\,dV=\int u{}\,\MQ_t\psi\,dV\).
For nonnegative $\psi$, divide by \({}t^2\), use \eqref{quadraticAveragingError}
and let \(t\to0\) to obtain the distributional bound.
\end{proof}

\begin{prop}[Interior Sobolev regularity]\label{prop:interiorSobolevRegularity}
Under the hypotheses of Theorem~\ref{thm:hessian},
\begin{equation}\label{W2pEstimate}
 u\in W^{2,p}_{\mathrm{loc}}(B_2)\qquad\text{for every }1<p<\infty.
\end{equation}
\end{prop}

\begin{proof}
Applying Proposition~\ref{prop:interiorLaplacian} to translates and dilates of
\(u\) gives ${}\Delta_e u\in L^\infty_{\mathrm{loc}}(B_2)$.
The interior \({}W^{2,p}\) estimate for the Poisson equation gives the conclusion.
\end{proof}

\subsection{The real Hessian}
\begin{proof}[Completion of the proof of Theorem~\ref{thm:hessian}]
For \({}x\in B_{3/2}\), apply Proposition~\ref{prop:interiorLaplacian} to
$v(y)=16u(x+y/4)$ on \({}B_2\). Then
${}[Dv]_{\alpha; B_2}\le4^\delta[Du]_{\alpha; B_2}$.
The function ${}g(x+\cdot/4)$ is bounded below by ${}g_*$, its Laplacian
is bounded below by $-L/16$, and its mean on ${}B_2$ is at most $4^{2n}$
times the mean of ${}g$ on $B_2$. Thus \(\Delta_eu\) is bounded on \(B_{3/2}\),
and \({}0\le v_t:=\MQ_t u/t^2\le C\) there for \({}0<t<1/8\).
Positivity and the Laplacian bound give ${}U=(u_{\ijbar})\in L^\infty(B_{3/2})$. On compact subballs, mollification gives bounded
\({}U_\varepsilon\to U\) in every finite \(L^p\), so the determinants converge
in ${}L^1$. Bedford-Taylor continuity under local uniform  approximation
identifies the limit with $g^n$. Thus $\det U=g^n$ almost everywhere.
Together  with $g\ge g_*$, the trace bound makes \(U\) uniformly elliptic
on \(B_{3/2}\).

Almost everywhere on \({}B_{3/2 }\), the arithmetic-geometric
mean inequality and \eqref{densityAverageLowerBound} imply
\begin{equation}\label{cofactorSubsolution}
 {}\mathcal Lv_t:=(\det U)\trace\big(U^{-1}((v_t)_{\ijbar})\big)
 \ge \frac{n g^n}{t^2}\left(\frac{\MG_t g}{g}-1\right)
 \ge -\frac{nL}{8(n+1)}g^{n-1}.
\end{equation}
This is a uniformly elliptic divergence-form inequality in the real
coordinates. Indeed the complex cofactor coefficients are divergence free
by commutation of mixed derivatives for smooth potentials and by
\(L^p_{\mathrm{loc}}\) convergence of  the mollified cofactors for \({}u\).
The corresponding real matrix is symmetric  and bounded. For fixed ${}t$,
${}\MG_t u\in C^2$ and Proposition~\ref{prop:interiorSobolevRegularity} give
\(v_t=(\MQ_t u)/t^2\in W^{2,p}_{\mathrm{loc}}\). Smooth approximation of $v_t$
therefore gives the weak divergence-form inequality.

Fix ${}x\in B_1$ and  write $M_x(s)=\sup_{B_s(x)}v_t$.
For \(0<s\le 1/8\), the nonnegative function \(M_x(2s)-v_t\) is a
supersolution on ${}B_{2s}(x)$ with bounded  source. The weak Harnack
inequality  with exponent one \cite[Theorem~8.18]{GT1983} gives
\begin{equation}\label{weakHarnackInequality}
 \frac{{}1}{|B_s|}\int_{B_s(x)}(M_x(2s)-v_t)
 \le C\big(M_x(2s)-M_x(s/2)+s^2\big).
\end{equation}
Here and below the constants depend only on $n$, ${}\alpha$,
${}|u|_{1,\alpha; B_2}$, $g_*$, ${}L$, and the mean of ${}g$ on ${}B_2$.
The exponent one is admissible in real dimension ${}2n$, since
$1<2n/(2n-2)$.

Let $\Gamma$ be the fundamental solution of \(\Delta_e\) in
\(\ER^{2n}\), and take a smooth cutoff \({}\eta\) supported in \({}B_{3/2}\)
and equal to one on \({}B_{5/4}\). For real indices \({}i,j\), the kernel
$D_{ij}\Gamma(y)$, away from zero, has size $C_n|y|^{-2n}$ and zero
integral on every centered annulus. On $s/2<|y|<s$ subtract \(M_x(2s)\)
from \(v_t(x-y)\). By \eqref{weakHarnackInequality},
\begin{equation}\label{annularKernelEstimate}
 \left|\int_{s/2<|y|<s} D_{ij}\Gamma(y)v_t(x-y)\,dy\right|
 \le C\big(M_x(2s)-M_x(s/2)+s^2\big).
\end{equation}
Sum over $s=2^{-k}/8$. The differences of suprema telescope, each
increment occurring  at most twice, and $\sum_k s^2<\infty$.
The part outside ${}|y|<1/8$ is bounded using the compact support of $\eta$.
It follows that
\begin{equation}\label{potentialHessianBound}
 {}\|D^2_{\ER}(\Gamma*(\eta v_t))\|_{L^\infty(B_1)}\le C .
\end{equation}
To justify the derivative, smooth \(\Gamma\) inside  a ball of dyadic
radius \(a\), keeping it unchanged outside and its second derivatives bounded
by \({}C_na^{-2n}\) inside. The inner contribution is bounded by
\({}C_n\|v_t\|_\infty\); passage to the distributional limit gives the estimate.

By the test-function limit above, \(v_t\to\Delta_eu/[8(n+1)]\)
in distributions. The same Hessian estimate therefore holds for
\({}\Gamma*(\eta\Delta_eu)\),  with the fixed dimensional factor.
The difference $u-\Gamma*(\eta\Delta_eu)$ is harmonic on
\(B_{5/4}\) and bounded there by the local bounds for \({}u,\Delta_eu\).
Interior harmonic estimates  bound its real Hessian on \(B_1\).
Adding the two Hessian bounds  gives ${}[Du]_{1; B_1}\le C$ by convexity
of ${}B_1$. Since ${}|u|_{1,\alpha; B_2}$ controls $\sup_{B_1}|u|$ and
$\sup_{{}B_1}|Du|$,  we obtain \eqref{realHessianEstimate}.
\end{proof}

\begin{proof}[Another Proof of Corollary~\ref{thm:main}]
Smooth positive $f$ satisfies the hypotheses of Theorem~\ref{thm:hessian}
on every ball whose closure is contained in ${}\Omega$. Proposition~\ref{prop:interiorSobolevRegularity}
therefore gives $u\in W^{2,p}_{\mathrm{loc}}$ for any 
\(p>n(n-1)\), so the interior regularity criterion of
B{\l}ocki and Dinew \cite[p.~412, (2)]{BD2011} gives smoothness.
\end{proof}

\appendix
\section{Gradient regularity from the Green averaging error}

\begin{lemma}\label{lem:greenErrorGradient}
Let \({}u\) be a continuous subharmonic function on \({}\Omega\subset\ER^m\),
\({}m\ge2\),  and let \({}B_R\Subset\Omega\). Suppose that \(0<a<1\) and
\[
 {}0\le \MQ_t u(x)=\MG_tu(x)-u(x)\le Kt^{1+a}
 \qquad (x\in B_R,\ 0<t<t_0),
\]
where \(t_0<\operatorname{dist}(B_R,\partial\Omega)\) and \(\MG_t\) is
defined by \eqref{greenAverageDefinition} using Euclidean ball averages.
Then  for every concentric ball \(B_r\), \(0<r<R\),
\[
 {}|u|_{1,a;B_r}\le C\bigl(|u|_{0;B_R}+K\bigr),
\]
where \({}C\) depends only on \(m,a,r,R,t_0\).
\end{lemma}

\begin{proof}
Let \({}\nu=\Delta_e u\), a  positive measure. For smooth \(u\), integration
of the normal derivative on spheres gives
\[
 \mathcal S_\sigma u(x)-u(x)
 =\frac{1}{|\mathbb S^{m-1}|}\int_0^\sigma
       \frac{\nu(B_\rho(x))}{\rho^{m-1}}\,d\rho .
\]
If \(s/2\le\sigma\le s<t_0\), integration over \(s/4<\rho<s/2\)
bounds the right side from below by \({}c_m s^{2-m}\nu(B_{s/4}(x))\).
This lower bound also holds for continuous subharmonic \({}u\). Indeed,
local mollifications converge uniformly on compact sets, and weak
convergence of their Laplacians gives lower semicontinuity of the mass
of the open ball \(B_{s/4}(x)\).

Since \({}\mathcal S_\sigma u(x)\ge u(x)\), integrating the lower bound
over \(s/2\le\sigma\le s\) in
\[
 {}\MB_su(x)=\frac{m}{s^m}\int_0^s \sigma^{m-1}\mathcal S_\sigma u(x)\,d\sigma
\]
gives \({}\MB_su(x)-u(x)\ge c_m s^{2-m}\nu(B_{s/4}(x))\).
Integrating over \({}t/2\le s\le t\) in \eqref{greenAverageDefinition} yields
\[
 {}c_m t^{2-m}\nu(B_{t/8}(x))\le\MQ_tu(x)\le Kt^{1+a}.
\]
Taking \({}t=8h\), we obtain
\[
 \nu(B_h(x))\le C K h^{m-1+a}
 \qquad (x\in B_R,\ 0<h<t_0/8).
\]

Choose \({}0\le\chi\le1\) smooth with compact support in \({}B_R\) and
equal to one  on \({}B_{(r+R)/2}\). Put \({}\eta=\chi\nu\).
Its total mass is \(\int u\Delta_e\chi\le |u|_{0;B_R}\int|\Delta_e\chi|\).
Every ball meeting \(\operatorname{supp}\chi{}\) lies in a ball of twice
its  radius centered in \(\operatorname{supp}\chi\). The preceding
estimates therefore give, with \({}M=K+|u|_{0;B_R}\),
\[
 \eta(B_h(x))\le CMh^{m-1+a}\qquad (x\in\ER^m,\ h>0).
\]

Let \(P=\Gamma*\eta\), where \({}\Delta_e\Gamma=\delta_0\).
Its distributional gradient is \(D\Gamma*\eta\). Since
\({}|D\Gamma(z)|\le C|z|^{1-m}\), summing over the annuli with radii
\(2^{-k-1}h<|x-y|\le2^{-k}h\) gives
\[
 \int_{{}B_h(x)} |D\Gamma(x-y)|\,d\eta(y)
 \le CM\sum_{k\ge0}(2^{-k}h)^a\le CMh^a .
\]
For \(\ell=|x-x'|>0\), this bounds the contribution of
\(B_{2\ell}(x)\cup B_{2\ell}(x')\) to \(|DP(x)-DP(x')|\)
by \(CM\ell^a\). On the  complement, the bound
\({}|D^2\Gamma(z)|\le C|z|^{-m}\) gives
\[
 |D\Gamma(x-y)-D\Gamma(x'-y)|\le C\ell |x-y|^{-m}.
\]
Summing over the annuli \({}2^k\ell\le |x-y|<2^{k+1}\ell\), \(k\ge1\),
bounds  the remaining contribution by
\(CM\ell\sum_{k\ge1}(2^k\ell)^{a-1}\le CM\ell^a\).
Thus \(DP\) has an \(a\)-H\"older continuous representative. Its supremum
on \(B_R\) is bounded by \({}CM\) by taking \({}h=2R\) in the near-field estimate.

Finally, \({}\Delta_e(u-P)=0\) on \({}B_{(r+R)/2}\).
Since \({}\Gamma\) is locally integrable and \({}\eta\) has mass at most \({}CM\),
Fubini's theorem gives \(\|P\|_{L^1(B_R)}\le CM\).
The interior derivative estimates for the harmonic function \({}u-P\)
therefore bound \({}|u-P|_{1,a;B_r}\) by \({}CM\).
Together with the estimates for \(DP\) and the bound for \(u\), this proves
the assertion.
\end{proof}

\section*{Acknowledgment}
AI Declaration. The authors proposed the initial idea for this paper. Through discussions with ChatGPT (GPT-6.2 Astra), they developed the strategy to control the Green averaging error  by computing at the contact point. The authors carefully checked  the entire manuscript and take full responsibility for its contents.

\noindent Jie Deng, Peking University, Beijing, China\\
\texttt{dj0401@stu.pku.edu.cn}\\

\noindent Jingchen Hu, Sichuan University, Chengdu, China\\
\texttt{jingchenhu@scu.edu.cn}\\

\noindent Liding Huang, Xiamen University, Xiamen, China\\
\texttt{huangliding@xmu.edu.cn}\\

\noindent Jiawei Liu, Nanjing University of Science \& Technology, Nanjing, China\\
\texttt{jiawei.liu@njust.edu.cn}\\

\noindent Li Sheng, Sichuan University, Chengdu, China\\
\texttt{lsheng@scu.edu.cn}\\

\noindent Xi Zhang, Nanjing University of Science \& Technology, Nanjing, China\\
\texttt{mathzx@njust.edu.cn}

\end{document}